\documentclass[letterpaper, 10pt, conference]{ieeeconf}      

\IEEEoverridecommandlockouts                              
\usepackage{graphics} 
\usepackage{epsfig} 
\usepackage{amsmath} 
\usepackage{amssymb}  

\let\bbbl\biggl

\let\bbbr\biggr

\newcommand{\norm}[1]{\lVert{#1}\rVert}

\newcommand{\R}{\mathbb{R}}

\newcommand{\half}{\frac{1}{2}}

\newcommand{\ip}[2]{\left\langle #1, #2 \right\rangle}

\newcommand{\bmat}[1]{\begin{bmatrix}#1\end{bmatrix}}

\newcommand{\mcl}[1]{\mathcal{ #1}}
\newcommand{\mbf}[1]{\mathbf{ #1}}

\newtheorem{theorem}{Theorem}

\newtheorem{lem}[theorem]{Lemma}

\title{\LARGE \bf
A New State-Space Representation for Coupled PDEs and Scalable Lyapunov Stability Analysis in the SOS Framework
}

\author{Matthew M. Peet
\thanks{This work was supported by grant NSF \# CMMI-1301851 and ONR \#N000014-17-1-2117}
\thanks{M. Peet is with School of Matter, Transport and Energy,
        Arizona State University, 501 Tyler Mall, Tempe, AZ, USA
        {\tt\small mpeet@asu.edu}}%
}

\begin{document}

\maketitle
\thispagestyle{empty}
\pagestyle{empty}

\begin{abstract}
We present a framework for stability analysis of systems of coupled linear Partial-Differential Equations (PDEs). The class of PDE systems considered in this paper includes parabolic, elliptic and hyperbolic systems with Dirichelet, Neuman and mixed boundary conditions. The results in this paper apply to systems with a single spatial variable and assume existence and continuity of solutions except in such cases when existence and continuity can be inferred from existence of a Lyapunov function. Our approach is based on a new concept of state for PDE systems which allows us to express the derivative of the Lyapunov function as a Linear Operator Inequality directly on $L_2$ and allows for any type of suitably well-posed boundary conditions. This approach obviates the need for integration by parts, spacing functions or similar mathematical encumbrances. The resulting algorithms are implemented in Matlab, tested on several motivating examples, and the codes have been posted online. Numerical testing indicates the approach has little or no conservatism for a large class of systems and can analyze systems of up to 20 coupled PDEs.  \vspace{-1mm}
\end{abstract}

\section{Introduction}\vspace{-1mm}

Partial Differential Equations (PDEs) are used to model systems where the state depends continuously on both time and secondary independent variables. Common examples of such secondary dependence include space, as in, e.g. rigid bodies (Bernoulli-Euler beams) and fluid flow (Navier-Stokes), or maturation, as in, e.g. cell populations and predator-prey dynamics.

The most common method for stability analysis of PDEs is to project the state onto a finite-dimensional vector space using, e.g. ~\cite{marion_1989,ravindran_2000,rowley_2005} and to use the existing extensive literature on control of ODEs to test stability and design controllers for the resulting finite-dimensional system. However, such discretization approaches are often prone to instability and numerical ill-conditioning. Attempts to develop a rigorous state-space theory for PDEs without discretization includes the significant literature on Semigroup theory~\cite{lasiecka_book,curtain_book,bensoussan_book}. Perhaps the most well-known method for stabilization of PDEs without discretization is the backstepping approach to controller synthesis~\cite{smyshlyaev_2005} (See the 2-state example in~\cite{aamo_2013}). Unfortunately, however, backstepping cannot currently be used for direct construction of Lyapunov functions for the purpose of stability analysis. Additional work on the use of computational methods and LMIs for computing Lyapunov functions for PDEs can be found in the work of~\cite{fridman_2009,fridman_2016,solomon_2015}. Other examples of LMI methods for stability analysis of PDEs include~\cite{gaye_2013}.

Recently, Sum-of-Squares (SOS) optimization methods have been applied to the problem of finding Lyapunov functions which prove stability of vector-valued PDEs. Examples of this work from our lab can be found in~\cite{papachristodoulou_2006CDC,gahlawat_2017TAC,gahlawat_2016ACC,gahlawat_2015CDC} and work from our colleagues can be found in~\cite{papachristodoulou_2006CDC,ahmadi_2016,valmorbida_2014,valmorbida_2016}. While these previous works have proven remarkably effective, they suffered from high computational complexity and the lack of a unifying framework - deficiencies which limit the practical impact and scalability of these results. The goal of this paper is to provide such a unifying framework and significantly reduce computational complexity by re-evaluating the state-space framework on which these earlier works were based.

Specifically, in this paper, we consider the problem of stability analysis of multiple coupled linear PDEs in a single spatial variable. We write these systems in the universal form
\[
\dot{x}(s,t) = A_0(s)x(s,t)+A_1(s)u_s(s,t)+A_2(s)x_{ss}(s,t)
\]
where $\mbf x$ is a \textbf{vector}-valued function $x:[a,b] \times \R^+ \rightarrow \R^n$ and with boundary constraints of the form
\[
B \bmat{ x(a,t) & x(b,t) & x_s(a,t) & x_s(b,t)}^T=0
\]
where $B$ is of row rank $2n$. These types of systems arise when there are multiple interacting spatially-distributed states and include wave equations, beam equations, et c.

The main technical result of this paper is to show that if $\mbf x$ satisfies the boundary conditions and is suitably differentiable, then we have the following identities
\begin{align*}
x(s)&= \int_a^b B_a(s,\eta) x_{ss}(\eta)d \eta+\int_a^s (s-\eta)  x_{ss}(\eta)d \eta\\
x_s(s)&= \int_a^b B_b(\eta) x_{ss}(\eta)d \eta+\int_a^s x_{ss}(\eta)d \eta,
\end{align*}
where the matrix-valued functions $B_a$ and $B_b$ are uniquely determined by the matrix $B$ and where $\mbf x_{ss}\in L_2[a,b]$ need not satisfy any constraints in order to define a solution. This identity implies that for any $\mbf x_{ss}$, the initial value problem is well-defined - implying that this is a boundary-condition independent representation of the state of the system.

We then use these identities to show that any Lyapunov function of the form
\begin{align*}
V(x)&=\int_a^b x(s)^T\bbbl( M(s)x(s) ds +\int_a^sN_1(s,\theta)x(\theta)d \theta\\[-2mm]
&\hspace{3.5cm}+\int_s^bN_2(s,\theta)x(\theta)d \theta\bbbr)ds
\end{align*}
may be equivalently represented as
\begin{align*}
&V(x)=\int_a^b x_{ss}(s)^T\bbbl( \int_a^s R_1(s,\theta)x_{ss}(\theta)d \theta\\[-3mm]
&\hspace{3.5cm}+\int_s^bR_2(s,\theta)x_{ss}(\theta)d \theta\bbbr)ds
\end{align*}
for some $R_1$, $R_2$ and furthermore, the derivative of this functional, $\dot V$, may likewise be represented in the same form. We note that the structure of these quadratic Lyapunov functions are implied by the closed-loop stability conditions established via the backstepping transformation, as shown in~\cite{gahlawat_2017TAC}. Furthermore, these results imply that the problem of computing stability of linear PDEs is equivalent to the problem of determining positivity of Lyapunov functions of this form for arbitrary functions $\mbf x_{ss}\in L_2$.

In the remainder of this paper, we will establish the results listed above, provide a computational framework for enforcing positivity of Lyapunov functions of this form, and show that the results are non-conservative and scalable through the use of numerical examples. Note that the identities listed can also be extended to third and fourth-order spatial derivatives, if required.

\section{Notation}
In this paper, we define $L_2^n[X]$ to be space of $\R^n$-valued Lesbegue integrable functions defined on $X$ and equipped with the standard inner product.
We use $W^{k,p}[X]$ to denote the Sobolev subspace of $L_p[X]$ defined as $\{u \in L_p[X]\, :\, \frac{\partial^q}{\partial x^q} u \in L_{p} \text{ for all } q \le k \}$.

\section{Preliminaries}\label{sec:prelims}
In this paper we consider stability of solutions $\mbf x:[a,b] \times \R^+ \rightarrow \R^n$ of PDEs of the form
\begin{equation}
x_t(s,t) = A_0(s)x(s,t)+A_1(s)x_s(s,t)+A_2(s)x_{ss}(s,t)\label{eqn:PDE_1}
\end{equation}
with boundary constraints of the form
\begin{equation}
B \bmat{ x(a,t)^T & x(b,t)^T & x_s(a,t)^T & x_s(b,t)^T}^T=0.\label{eqn:PDE_2}
\end{equation}
These boundary conditions can be used to represent Dirichelet, Neumann, Robin, et c., with the only restriction that the row rank of $B$ need be $2n$. In the semigroup framework, this translates to $\dot{\mbf x}=\mathcal{A}\mbf x$ with generator
\[
\mathcal{A}:=A_0(s)+A_1(s)\partial_s+A_2(s)\partial_{ss}
\]
and domain
\begin{align*}
&D_A:=\\
&\{\mbf x\in W^{2,2}\,:\, B \bmat{ \mbf x(a)^T & \hspace{-1mm}\mbf x(b)^T & \hspace{-1mm}\mbf x_s(a)^T & \hspace{-1mm}\mbf x_s(b)^T}^T=0 \}
\end{align*}

\section{Lyapunov Stability}\label{sec:Lyap}
It seems that existence of a Lyapunov function does not guarantee existence and continuity of solutions for PDEs except in certain very limited special cases. Therefore, we must assume these properties hold and we give mathematical rigour to this assumption by assuming the existence of a ``Semi-continuous semigroup'', $S(t):X\rightarrow X$ with domain $D_A\subset X$ so that $S(\tau)\mbf x(\cdot,t)=\mbf x(\cdot,t+\tau)$ for any solution to Eqns.~\eqref{eqn:PDE_1} and~\eqref{eqn:PDE_2}. See~\cite{lasiecka_book}. The following is from~\cite{curtain_book}.

\begin{theorem}\label{thm:stability}
Suppose that $\mathcal{A}$ generates a strongly-continuous semigroup on $X$ with domain $D_A$ and there exists $\alpha,\beta,\gamma>0$ and $\mathcal{P}:X \rightarrow X$ such that $\alpha \norm{\mbf x}_X \le \ip{\mbf x}{\mathcal{P}\mbf x}_X\le \beta \norm{\mbf x}_X$ and
\[
\ip{\mbf x}{\mathcal{PA}\mbf x}_X+\ip{\mathcal{A}\mbf x}{\mathcal{P}\mbf x}_X\le -\gamma \norm{\mbf x}_X
\]
for all $\mbf x \in D_A$. Then the system defined by Eqns~\eqref{eqn:PDE_1} and~\eqref{eqn:PDE_2} is exponentially stable in $\norm{\cdot}_X$.
\end{theorem}
In this paper, we show how these conditions may be enforced when $\mcl A$ and $D_{\mcl A}$ are as defined in Section~\ref{sec:prelims}. In this case $\mcl A$ is a differential operator. We will show in Sections~\ref{sec:identity} and~\ref{sec:reform} that the stability conditions in Theorem~\ref{thm:stability} can be reformulated with on $D_{A}=L_2$ and in Sections~\ref{sec:LOI} and~\ref{sec:LOI2} we will show that these conditions can be enforced using LMIs based on an SOS-style approach.

\section{Fundamental Identities}\label{sec:identity}

In this section, we show that if
\[
B \bmat{ x(a)^T & x(b)^T & x_s(a)^T & x_s(b)^T}^T=0
\]
where $B$ is of row rank $2n$, then the following identities hold
\begin{align*}
x(s)&= \int_a^b B_a(s,\eta) x_{ss}(\eta)d \eta+\int_a^s (s-\eta)  x_{ss}(\eta)d \eta\\
x_s(s)&= \int_a^b B_b(\eta) x_{ss}(\eta)d \eta+\int_a^s x_{ss}(\eta)d \eta,
\end{align*}
where $B_a$ and $B_b$ are uniquely determined by the matrix $B$.

First, we establish the auxiliary identities:

\begin{lem}\label{lem:identity1}Suppose that $x$ is twice continuously differentiable. Then
\begin{align*}
x_s(s)&=x_s(a)+\int_a^s x_{ss}(\eta)d \eta\\
x(s)&=x(a)+x_s(a)(s-a)+\int_a^s (s-\eta)  x_{ss}(\eta)d \eta
\end{align*}
\end{lem}
\begin{proof}
  The first identity is the fundamental theorem of calculus. The second identity is a repeated application of the fundamental theorem of calculus, combined with a change of variables. That is,
\begin{align*}
x(s)&=x(a)+\int_a^s  x_{s}(\eta)d \eta\\
&=x(a)+\int_a^s x_s(a)ds+ \int_a^s  \int_a^\eta x_{ss}(\zeta)d \zeta d \eta
\end{align*}
Examining the 3rd term, where $I(s)$ is the indicator function,

\begin{align*}
&\int_a^s  \int_a^\eta x_{ss}(\zeta)d \zeta d \eta=\int_a^b  \int_a^b I(s-\eta)I(\eta-\zeta) x_{ss}(\zeta)d \zeta d \eta\\
&=\int_a^b  \left(\int_a^b I(s-\eta)I(\eta-\zeta) d \eta\right) x_{ss}(\zeta) d \zeta\\
&=\int_a^b  I(s-\zeta)\left(\int_s^\zeta  d \eta\right) x_{ss}(\zeta) d \zeta= \int_a^s  \left(s-\zeta\right) x_{ss}(\zeta) d \zeta
\end{align*}
which is the desired result
\end{proof}

As an obvious corollary, we have
\begin{align*}
x_s(b)&=x_s(a)+\int_a^b x_{ss}(\eta)d \eta\\
x(b)&=x(a)+x_s(a)(b-a)+\int_a^b (b-\eta)  x_{ss}(\eta)d \eta \\
\end{align*}


The implication is that any boundary value can be expressed using two other boundary identities. We can now generalize this to the main result.

\begin{theorem}\label{thm:identity2} Suppose $\mbf x\in W^{3,2}[a,b]$ and
\[
B \bmat{ x(a)^T & x(b)^T & x_s(a)^T & x_s(b)^T}^T=0
\]
where $B$ has row rank $2n$, then
\begin{align*}
x(s)&= \int_a^b B_a(s,\eta) x_{ss}(\eta)d \eta+\int_a^s (s-\eta)  x_{ss}(\eta)d \eta\\
x_s(s)&= \int_a^b B_b(\eta) x_{ss}(\eta)d \eta+\int_a^s x_{ss}(\eta)d \eta,
\end{align*}
where
\begin{align*}
&B_a(s,\eta)=B_4(s)(b-\eta)+B_5(s),\\
&B_b(\eta)=B_6(b-\eta) +  B_7\\
&\bmat{B_6 & B_7}  =\bmat{0& I} B_3,\\
& \bmat{B_4(s) & B_5(s)}  =\bmat{I&(s-a)I} B_3\\
    B_3&=B_2^{-1}B \bmat{0 & 0\\ I & 0\\ 0 & 0\\ 0& I},\qquad B_2=B\bmat{I &0\\I&(b-a)I\\0&I\\0 &I}
\end{align*}
\end{theorem}
\begin{proof}
Using Lemma~\ref{lem:identity1}, we can express all boundary terms using $x(a)$, $x_s(a)$, and $x_{ss}(s)$.
\[
\bmat{x(a)\\x(b)\\x_s(a)\\x_s(b)}= \bmat{I &0\\I&(b-a)I\\0&I\\0 &I}\bmat{x(a)\\x_s(a)}+  \bmat{ 0\\ \int_a^b (b-\eta) x_{ss}(\eta)d \eta\\0\\ \int_a^b x_{ss}(\eta)d \eta}
\]
Hence
\begin{align*}
&B\bmat{x(a)\\x(b)\\x_s(a)\\x_s(b)}= \underbrace{B\bmat{I &0\\I&(b-a)I\\0&I\\0 &I}}_{B_2}\bmat{x(a)\\x_s(a)}\\
&\qquad \qquad \qquad \qquad  +  B \bmat{ 0\\ \int_a^b (b-\eta) x_{ss}(\eta)d \eta\\0\\ \int_a^b x_{ss}(\eta)d \eta}\\
&= B_2\bmat{x(a)\\x_s(a)}+  B \bmat{0 & 0\\ I & 0\\ 0 & 0\\ 0& I}\bmat{  \int_a^b (b-\eta) x_{ss}(\eta)d \eta\\ \int_a^b x_{ss}(\eta)d \eta}=0
\end{align*}
Since $B$ has $2n$ row rank, $B_2$ is invertible and hence we have
\[
\bmat{x(a)\\x_s(a)}=  \underbrace{-B_2^{-1}B \bmat{0 & 0\\ I & 0\\ 0 & 0\\ 0& I}}_{B_3}\bmat{  \int_a^b (b-\eta) x_{ss}(\eta)d \eta\\ \int_a^b x_{ss}(\eta)d \eta}
\]

Now, using Lemma~\ref{lem:identity1},{\small
\begin{align*}
&x(s)=x(a)+x_s(a)(s-a)+\int_a^s (s-\eta) x_{ss}(\eta)d \eta \\
&=\bmat{I&(s-a)I} \bmat{x(a)\\x_s(a)} +\int_a^s (s-\eta)  x_{ss}(\eta)d \eta \\
&=\underbrace{\bmat{I&(s-a)I} B_3}_{\bmat{B_4(s) & B_5(s)}}\bmat{  \int\limits_a^b (b-\eta) x_{ss}(\eta)d \eta\\ \int_a^b x_{ss}(\eta)d \eta} +\hspace{-1mm}\int\limits_a^s (s-\eta)  x_{ss}(\eta)d \eta \\
&= \int_a^b B_4(s) (b-\eta) x_{ss}(\eta)d \eta\\[-3mm]
& \qquad \qquad \qquad \qquad +  B_5(s)\int_a^b x_{ss}(\eta)d \eta +\int_a^s (s-\eta)  x_{ss}(\eta)d \eta \\
&= \int_a^b \left(B_4(s)(b-\eta)+B_5(s)\right) x_{ss}(\eta)d \eta+\int_a^s (s-\eta)  x_{ss}(\eta)d \eta \\
&= \int_a^b B_a(s,\eta) x_{ss}(\eta)d \eta+\int_a^s (s-\eta)  x_{ss}(\eta)d \eta.
\end{align*}}
Likewise, we have
\begin{align*}
x_s(s)&=x_s(a)+\int_a^s x_{ss}(\eta)d \eta\\
&=\bmat{0&I}\bmat{x(a)\\x_s(a)}+\int_a^s x_{ss}(\eta)d \eta\\
&=\underbrace{\bmat{0&I}B_3}_{\bmat{B_6 & B_7}} \bmat{  \int_a^b (b-\eta) x_{ss}(\eta)d \eta\\ \int_a^b x_{ss}(\eta)d \eta}+\int_a^s x_{ss}(\eta)d \eta\\
&= \int_a^b \left(B_6(b-\eta) +  B_7\right)x_{ss}(\eta)d \eta+\int_a^s x_{ss}(\eta)d \eta\\
&= \int_a^b B_b(\eta) x_{ss}(\eta)d \eta+\int_a^s x_{ss}(\eta)d \eta.
\end{align*}
\end{proof}
\section{Reformulation of the Lyapunov Function}\label{sec:reform}

If we denote the class of operators $\mcl P_{\{M,N_1,N_2\}}:L_2^n \rightarrow L_2^n$ by
\begin{align*}
&\left(\mcl P_{\{M,N_1,N_2\}}\mbf x\right)(s)\\
&= M(s) x(s) ds +\int\limits_a^s N_1(s,\theta)x(\theta)d \theta+\int\limits_s^bN_2(s,\theta)x(\theta)d \theta,
\end{align*}
then we may compactly represent our Lyapunov candidate form as
\[
V(\mbf x)=\ip{\mbf x}{\mcl P_{\{M,N_1,N_2\}}\mbf x}_{L_2}
\]
The derivative of the Lyapunov candidate may then be likewise compactly represented (w/ slight abuse of notation) as
\begin{align*}
&\half \dot V(\mbf x)=\ip{\mbf x}{\mcl P_{\{MA_0,N_1A_0,N_2A_0\}}\mbf x}\\
&+\ip{\mbf x}{\mcl P_{\{MA_1,N_1A_1,N_2A_1\}}\mbf x_s}+\ip{\mbf x}{\mcl P_{\{MA_2,N_1A_2,N_2A_2\}}\mbf x_{ss}}
\end{align*}
The challenge, then, is to show that each of these terms may, in turn, be represented in the form
\[
\ip{\mbf x_{ss}}{\mcl P_{\{0,R_1,R_2\}}\mbf x_{ss}}_{L_2}
\]
through repeated use of the identities
\begin{align*}
x(s)&= \int_a^b B_a(s,\eta) x_{ss}(\eta)d \eta+\int_a^s (s-\eta)  x_{ss}(\eta)d \eta\\
x_s(s)&= \int_a^b B_b(\eta) x_{ss}(\eta)d \eta+\int_a^s x_{ss}(\eta)d \eta.
\end{align*}

For convenience, we leave off the $A_i$ terms and address each inner product separately. Let use also define the following functions which are common to all three results.{\small
\begin{align*}
Y_1(s,\eta)&=B_a(\eta,s)^T M(\eta) + \int_\eta^b  B_a(\theta,s)^T N_1(\theta,\eta)d\theta\\
&\hspace{3.5cm}+\int_a^\eta B_8(\theta,s)^T N_2(\theta,\eta)d\theta\\
Y_2(\zeta)&=M(\zeta)  + \int_a^\zeta    N_1(\zeta,\theta)   d \theta + \int_\zeta^b  N_2(\zeta,\theta)   d \theta\\
Y_3(\zeta,\eta)&=  M(\zeta)B_a(\zeta,\eta) + \int_a^\zeta   N_1(\zeta,\theta)  B_a(\theta,\eta)  d \theta \\
&\hspace{3.5cm}+  \int_\zeta^b N_2(\zeta,\theta) B_a(\theta,\eta)  d \theta
\end{align*}}
Note that these functions are defined in terms of $M$, $N_1$, and $N_2$ and hence will vary if these terms are defined differently for Lemmas~\ref{lem:reform1},~\ref{lem:reform2}, and~\ref{lem:reform3}.

\begin{lem}\label{lem:reform1}Suppose $\mbf x$ satisfies the conditions of Thm.~\ref{thm:identity2}. Then\vspace{-3mm}
\[
 \ip{\mbf x}{\mcl P_{\{M,N_1,N_2\}}\mbf x_{ss}}= \ip{\mbf x_{ss}}{\mcl P_{\{0,R_1,R_2\}}\mbf x_{ss}}\vspace{-2mm}
\]
where\vspace{-2mm}
\begin{align*}
R_1(s,\theta)&=E_1(s,\theta)+E_3(s,\theta),\\
R_2(s,\theta)&=E_2(s,\theta) +E_3(s,\theta)
\end{align*}
\begin{align*}
&E_1(s,\theta)=\int_s^b (\eta-s)  N_1(\eta,\theta)d\eta\\
&E_2(s,\theta)=(\theta-s)M(\theta)+\int_\theta^b (\eta-s) N_1(\eta,\theta)d\eta\\
&\hspace{4.5cm} +\int_s^\theta(\eta-s) N_2(\eta,\theta)d\eta \\
&E_3(s,\theta) =Y_1(s,\theta)\\
\end{align*}
\end{lem}
\begin{proof} The proofs of these Lemmas cannot be included in conference format due to length constraints. Therefore, for the proof of these lemmas, we refer to an Arxiv Appendix, available online at~\cite{peet_arxiv}
\end{proof}

\noindent \textbf{Notation:} For convenience, we say
\[
(R_1,R_2)=\mcl L_1(M,N_1,N_2)
\]
if $R_1$, $R_2$, $M$, $N_1$, and $N_2$ satisfy the conditions of Lemma~\ref{lem:reform1}.

\begin{lem}\label{lem:reform2}Suppose $\mbf x$ satisfies the conditions of Thm.~\ref{thm:identity2}. Then\vspace{-3mm}
\[
 \ip{\mbf x}{\mcl P_{\{M,N_1,N_2\}}\mbf x_s}= \ip{\mbf x_{ss}}{\mcl P_{\{0,Q_1,Q_2\}}\mbf x_{ss}}\vspace{-2mm}
\]
where\vspace{-2mm}
\begin{align*}
Q_1(s,\theta)&=F_1(s,\theta)+F_3(s,\theta)\\
Q_2(s,\theta)&=F_2(s,\theta) +F_3(s,\theta)
\end{align*}
\begin{align*}
&F_1(s,\theta)= \int_s^b\left(  (\eta-s)F_4(\theta,\eta)+ F_5(s,\eta) \right)d \eta\\
&F_2(s,\theta)=\int_\theta^b \left( (\eta-s) F_4(\theta,\eta) + F_5(s,\eta)\right)d \eta\\
&F_3(s,\eta) =\int_a^b B_a(\zeta,s)^T Y_2(\zeta) B_b(\eta) d \zeta +\int_\eta^b Y_1(s,\zeta)d \zeta\\
&\hspace{4.5cm} +\int_s^b (\zeta-s)Y_2(\zeta) d \zeta B_b(\eta)\\
&F_4(\theta,\eta)= M(\eta)  +  \int_\theta^\eta   N_1(\eta,\zeta)   d \zeta\\
&F_5(s,\eta)= \int_s^\eta  (\zeta - s)  N_2(\zeta,\eta)    d \zeta.
  \end{align*}
\end{lem}

\noindent \textbf{Notation:} For convenience, we say
\[
(Q_1,Q_2)=\mcl L_2(M,N_1,N_2)
\]
if $Q_1$, $Q_2$, $M$, $N_1$, and $N_2$ satisfy the conditions of Lemma~\ref{lem:reform2}.

\begin{lem}\label{lem:reform3}Suppose $\mbf x$ satisfies the conditions of Thm.~\ref{thm:identity2}. Then\vspace{-3mm}
\[
 \ip{\mbf x}{\mcl P_{\{M,N_1,N_2\}}\mbf x}= \ip{\mbf x_{ss}}{\mcl P_{\{0,T_1,T_2\}}\mbf x_{ss}}\vspace{-2mm}
\]
where
\begin{align*}
T_1(s,\theta)&=G_1(s,\theta)+G_3(s,\theta)\\
T_2(s,\theta)&=G_2(s,\theta) +G_3(s,\theta)
\end{align*}

\begin{align*}
&G_1(s,\theta)=\int_s^b \left((\eta-s) G_4(\theta,\eta) +G_5(s,\theta,\eta) \right) d \eta \\
&G_2(s,\theta)=\int_\theta^b \left((\eta-s) G_4(\theta,\eta)+ G_5(s,\theta,\eta) \right) d \eta\\
&G_3(s,\theta)=\int_a^b B_a(\eta,s)^T Y_3(\eta,\theta)d\eta\\
&\hspace{1.3cm}+\int_\theta^b (\eta-\theta) Y_1(s,\eta) d \eta+\int_s^b (\eta-s) Y_3(\eta,\theta) d\eta\\
&G_4(\theta,\eta)= (\eta-\theta) M(\eta)  +  \int_\theta^\eta  (\zeta-\theta) N_1(\eta,\zeta)  d \zeta \\
&G_5(s,\theta,\eta)=\int_s^\eta  (\zeta-s) (\eta-\theta)  N_2(\zeta,\eta)   d \zeta.
\end{align*}
\end{lem}

\noindent \textbf{Notation:} For convenience, we say
\[
(T_1,T_2)=\mcl L_3(M,N_1,N_2)
\]
if $T_1$, $T_2$, $M$, $N_1$, and $N_2$ satisfy the conditions of Lemma~\ref{lem:reform3}.

Note that the operators obtained here are not necessarily symmetric. However, we may construct an equivalent symmetric representation as $\mcl P + \mcl P^*$ using
\[
\mcl P_{\{0,Q_2(\theta,s)^T,Q_1(\theta,s)^T\}}^*=\mcl P_{\{0,Q_1(s,\theta),Q_2(s,\theta)\}}.
\]
That is, in the symmetric representation, $Q_1(s,\theta)= Q_2(\theta,s)^T$.

\section{Positivity of Operators}\label{sec:LOI}
Now that we have shown how to represent our Lyapunov stability conditions as positivity of operators of the form $\mcl P_{\{0,N_1,N_2\}}$, we must show how to use LMIs to enforce positivity of these operators when $N_1$ and $N_2$ are polynomials. This is a slight generalization of the result in~\cite{peet_2014ACC}.

\begin{theorem}\label{thm:LOI}
For any square-integrable functions $Z(s)$ and $Z(s,\theta)$, if $g(s)\ge 0$ for all $s \in [a,b]$ and{\small
\begin{align*}
&M(s)=g(s)Z(s)^T P_{11}Z(s)\\
&N_1(s,\theta)=g(s) Z(s)^T P_{12} Z(s,\theta) + g(\theta)Z(\theta,s)^T P_{31} Z(\theta)\\
&  + \hspace{-1mm}\int_{0}^\theta\hspace{-1.5mm} g(\nu) Z(\nu,s)^T P_{33} Z(\nu,\theta) d\nu   +\hspace{-1mm}\int_\theta^s \hspace{-1.5mm}g(\nu)Z(\nu,s)^T P_{32} Z(\nu,\theta) d\nu\\
&       +\int_s^{L}g(\nu) Z(\nu,s)^T P_{22} Z(\nu,\theta) d\nu\\
&N_2(s,\theta)= g(s) Z(s)^T P_{13} Z(s,\theta) + g(\theta)Z(\theta,s)^T P_{21} Z(\theta)\\
&  + \int_{0}^s\hspace{-1.5mm}g(\nu) Z(\nu,s)^T P_{33} Z(\nu,\theta) d\nu
       +\int_s^\theta \hspace{-1.5mm}g(\nu)Z(\nu,s)^T P_{23} Z(\nu,\theta) d\nu\\
&       +\int_\theta^{L} g(\nu) Z(\nu,s)^T P_{22} Z(\nu,\theta) d\nu,
\end{align*}}
where
\[
P = \bmat{ P_{11} & P_{12}& P_{13}\\
      P_{21} & P_{22}& P_{23}\\
      P_{31} & P_{32}& P_{33}}\ge 0,\vspace{-2mm}
\]
then $\ip{\mbf x}{\mcl P_{\{M,N_1,N_2\}}\mbf x}_{L_2}\ge 0$ for all $\mbf x \in L_2[a,b]$.
\end{theorem}
\begin{proof}
Define the operator
\[
\left(\mcl Z \mbf x\right)(s) = \bmat{\sqrt{g(s)}Z(s)\mbf x(s)\\
\int_0^s \sqrt{g(s)}Z(s,\theta)\mbf x(\theta)d\theta\\
\int_s^b \sqrt{g(s)} Z(s,\theta)\mbf x(\theta)d\theta}.
\]
Then\vspace{-1mm}
 \[
\ip{\mbf x}{\mcl P_{\{M,N_1,N_2\}} \mbf x}=\ip{\mcl Z \mbf x}{P \mcl Z \mbf x}=\ip{P^{\half} \mcl Z \mbf x}{P^{\half} \mcl Z \mbf x}\ge 0.
\]

\end{proof}

For convenience, we define the cone of such operators as
\begin{align}
\Phi:=&\{(M,N_1,N_2)\,:\, \text{$M$, $N_1$ and $N_2$ satisfy} \notag \\
&\qquad \qquad   \text{the conditions of Thm.~\ref{thm:LOI}.}\}\label{defn:Phi}
\end{align}
where the dimension of the matrices $M$, $N_1$ and $N_2$ should be clear from context.

\section{SOS Conditions for Stability}\label{sec:LOI2}

The stability conditions can now be written concisely using the definitions of $\Phi$, $\mcl L_1$, $\mcl L_2$, and $\mcl L_3$ as follows.
\begin{theorem}\label{thm:main}
Suppose there exist $\epsilon>0$,
\[
(M-\epsilon I,N_1,N_2) \in \Phi\]
and
\[
\left(0,-H_1(s,\theta)-H_2(\theta,s)^T,-H_2(s,\theta)-H_1(\theta,s)^T\right)\in \Phi
\]
where
\begin{align*}
&(H_1,H_2)=\mcl L_1(V_1,W_{11},W_{12})\\
&\hspace{2cm}+\mcl L_2(V_2,W_{21},W_{22})+\mcl L_3(V_3,W_{31},W_{32})
\end{align*}
\begin{align*}
&V_1(s)=M(s)A_0(s)+\epsilon I, \quad &&W_{11}(s,\theta)=N_1(s,\theta)A_0(\theta),\\
&\hspace{3cm} &&W_{12}(s,\theta)=N_2(s,\theta)A_0(\theta)\\
&V_2(s)=M(s)A_1(s),\quad &&W_{21}(s,\theta)=N_1(s,\theta)A_1(\theta),\\
&\hspace{3cm} &&W_{22}(s,\theta)=N_2(s,\theta)A_1(\theta))\\
&V_3(s)=M(s)A_2(s),\quad &&W_{31}(s,\theta)=N_1(s,\theta)A_2(\theta),\\
& && W_{32}(s,\theta)=N_2(s,\theta)A_2(\theta))
\end{align*}
Then any solution of Eqns.~\eqref{eqn:PDE_1} and~\eqref{eqn:PDE_2} is exponentially stable.
\end{theorem}
\begin{proof}
Let
\[
V(\mbf x)=\ip{\mbf x}{\mcl P_{\{M,N_1,N_2\} }\mbf x}\ge  \epsilon \norm{\mbf x}^2_{L_2}.
\]
Then
\begin{align*}
&\dot V(\mbf x)+2\epsilon \norm{\mbf x}^2_{L_2}=
2\ip{\mbf x}{P_{\{V_1,W_{11},W_{12}\}} \mbf x}\\
&+2\ip{\mbf x}{P_{\{V_2,W_{21},W_{22}\}} \mbf x_s}
+2\ip{\mbf x}{P_{\{V_3,W_{31},W_{32}\}} \mbf x_{ss}}\\
&=2\ip{\mbf x_{ss}}{\mcl P_{\{0,H_1,H_2\}} \mbf x_{ss}}\le 0
\end{align*}
Therefore, we have exponential stability from Thm.~\ref{thm:stability}.
\end{proof}

\section{Numerical Implementation and Analysis}\label{sec:examples}
In this section, we examine the accuracy and computational complexity of the proposed stability algorithm by applying the results to several well-studied problems. The algorithms are implemented using a Matlab toolbox which is an adaptation of SOSTOOLS~\cite{prajna_2002} and which can be found online at \texttt{http://control.asu.edu}. In all cases, the conditions of Theorems~\ref{thm:main} and~\ref{thm:LOI} are applied by choosing $Z$ to be a vector of monomial bases of degree $d$ and less and either $g(s)=1$ or $g(s)=(s-a)(b-s)$.

\noindent\textbf{Example 1:} We begin with several variations of the diffusion equation. The first is adapted from~\cite{valmorbida_2014}.

\[
\dot x(t,s)=\lambda x(t,s) + x_{ss}(t,s)\vspace{-2mm}
\]
where $x(0)=x(1)=0$ and which is known to be stable if and only if $\lambda <\pi^2=9.8696$. For $d=1$, the algorithm is able to prove stability for $\lambda=9.8696$ with a computation time of .54s.

\noindent\textbf{Example 2:} The second example from~\cite{valmorbida_2016} is the same, but changes the boundary conditions to $x(0)=0$ and $x_s(1)=0$ and is unstable for $\lambda>2.467$. For $d=1$, the algorithm is able to prove stability for $\lambda=2.467$ with identical computation time.

\noindent\textbf{Example 3:} The third example from~\cite{gahlawat_2017TAC} is not homogeneous
\begin{align*}
\dot x(t,s)=&(-.5 s^3+1.3 s^2-1.5 s+.7+\lambda) x(t,s)\\
& + (3s^2-2s)x_{s}(t,s) + (s^3-s^2+2)x_{ss}(t,s)\vspace{-2mm}
\end{align*}
where $x(0)=0$ and $x_s(1)=0$ and was estimated numerically to be unstable for $\lambda > 4.65$. For $d=1$, the algorithm is able to prove stability for $\lambda=4.65$ with similar computation time.

\noindent\textbf{Example 4:} In this example from~\cite{valmorbida_2014}, we have
\[
\dot x(t,s)=\bmat{1 & 1.5 \\ 5 & .2}x(t,s)+ R^{-1} x_{ss}(t,s)\vspace{-2mm}
\]
with $x(0)=0$ and $x_s(1)=0$. In this case, using $d=1$, we can prove stability for $R=2.93$ (improvement over $R=2.45$ in~\cite{valmorbida_2014}) with a computation time of $1.21s$.

\noindent\textbf{Example 5:} In this example from~\cite{valmorbida_2016}, we have
\[
\dot x(t,s)=\bmat{0 & 0 & 0\\ s & 0 & 0\\ s^2 & -s^3& 0}x(t,s)+ R^{-1} x_{ss}(t,s)\vspace{-2mm}
\]
with $x(0)=0$ and $x_s(1)=0$. In this case, using $d=1$, we prove stability for $R=21$ (and greater) with a computation time of $4.06s$.

\noindent\textbf{Example 6:} Next, we consider a damped wave equation $x_{tt}=x_{ss}-kx_t$ with $x(0)=x(1)=0$

\[
\dot x(t,s)=\bmat{0 & 1\\ 0 & -k}x(t,s) + \bmat{0 & 0\\1 & 0} x_{ss}(t,s)\vspace{-2mm}
\]
This is shown to be stable for $k=.1$ with a computation time of $1.54s$.

\noindent\textbf{Example 7:} Finally, we explore computational complexity using a simple $n$-dimensional diffusion equation

\[
\dot x(t,s)=x(t,s) +  x_{ss}(t,s)\vspace{-2mm}
\]
where $x(t,s) \in \R^n$. We then evaluate the computation time for different size problems, from $n=1$ to $n=20$.

\[
\hbox{\begin{tabular}{c|c|c|c|c}\label{tab:taumax}
$n$ & $1$  & $5$ & $10$ & $20$ \\
\hline
CPU sec & .54 & 37.4 & 745 & 31620  \\
\end{tabular}}
\]

\section{Conclusion} In this paper, we have shown that stability of a large class of PDE systems can be represented compactly in LMI form using a variation of Sum-of-Squares optimization. To achieve this result, we proposed that the state of a PDE of the form of Equation~\eqref{eqn:PDE_1} is actually $\mbf x_{ss}$ and that all Lyapunov stability conditions may be represented on this state. A SOS-style algorithm to test these Lyapunov conditions is proposed and numerical examples indicate no conservatism in the stability conditions to at least $5$ significant figures even for low polynomial degree. It is clear that these results can also be directly extended to: PDEs with uncertainty; $H_{\infty}$-gain analysis of PDEs;  $H_\infty$-optimal observer synthesis for PDEs; $H_\infty$-optimal control of PDEs. Nonlinear Stability analysis can likewise be considered. In addition, the identities proposed in Section~\ref{sec:identity} by be extended to multiple spatial dimensions. Some unanswered questions include how to repose several common stability problems in the proposed generalized framework. For example, the wave equation with $u_t(L)=-u_x(L)$ in its native form is not suitably well-posed as the $B$ matrix does not have row rank $2n$.

\bibliographystyle{IEEEtran}
\bibliography{peet_bib,delay,PDEs}
\end{document}